\documentclass[12pt]{amsart}
\usepackage{pstricks}
\usepackage{amsfonts}
\usepackage{amssymb}
\usepackage{amsthm}
\usepackage{latexsym,amsmath}
\usepackage{graphicx}
\usepackage{stmaryrd}
\usepackage{natbib}
\setcitestyle{numbers,aysep={,},yysep={;},notesep={, }}

\newtheorem{theorem}{Theorem}[section]
\theoremstyle{definition}
\newtheorem{Definition}[theorem]{Definition}

\newenvironment{theorem*}[1]{\medskip
                            \noindent
                            {\bf Theorem #1. }\ %
                            \begingroup \sl}
                            {\endgroup\medskip}

\title
 {Intrinsic justifications for large-cardinal axioms}

\author[$\mathrm{M^{\lowercase{c}}Callum}$ ]{\textbf{Rupert} $\mathbf{M^{\lowercase{c}}Callum}$ }

\begin{document}

\begin{abstract}

We shall defend three philosophical theses about the extent of intrinsic justification based on various technical results. We shall present a set of theorems which indicate intriguing structural similarities between a family of ``weak" reflection principles roughly at the level of those considered by W. Tait and P. Koellner and a family of ``strong" reflection principles roughly at the level of those of P. Welch and S. Roberts, which we claim to lend support to the view that the stronger reflection principles are intrinsically justified as well as the weaker ones. We consider connections with earlier work of V. Marshall.

\end{abstract}

\maketitle

\bigskip

\section*{Acknowledgements}

I am thankful to Neil Barton for helpful feedback on a draft version of this work.

\section{Introduction}

In speaking of justifications for candidates for new axioms for set theory to be added to the ZFC axioms, philosophers of set theory draw a distinction between ``intrinsic" and ``extrinsic" justifications. This distinction has its origins in remarks made by Kurt G\"odel in his famous essay about the continuum problem \cite{Godel1964}. There he wrote

\bigskip

\begin{quotation}
First of all, the axioms of set theory by no means form a system closed in itself, but, quite on the contrary, the very concept of set on which they are based suggests their extension by new axioms which assert the existence of still further iterations of the operation ``set of$\,$". These axioms can be formulated also as propositions asserting the existence of very great cardinal numbers (i.e., of sets having these cardinal numbers). The simplest of these strong ``axioms of infinity" asserts the existence of inaccessible numbers (in the weaker or stronger sense) $> \aleph_{0}$. The latter axiom, roughly speaking, means nothing else but that the totality of sets obtainable by use of the procedures of formation of sets expressed in the other axioms forms again a set (and, therefore, a new basis for further applications of these procedures). Other axioms of infinity have first been formulated by P. Mahlo. These axioms show clearly, not only that the axiomatic system of set theory as used today is incomplete, but also that it can be supplemented without arbitrariness by new axioms which only unfold the content of the concept of set explained above.'
\end{quotation}

\bigskip

\noindent Thus here he indicated that he took the view that inaccessible and Mahlo cardinals were intrinsically justified and `implied by our general conception of set'.

\bigskip

\noindent In a footnote to the 1964 version of the essay, which was inserted in September 1966, he wrote `... some propositions have been formulated which, if consistent, are extremely strong axioms of infinity of an entirely new kind... That these axioms are implied by the general concept of set in the same sense as Mahlo's has not been made clear yet.' Here he was referring to axioms asserting the existence of weakly compact, measurable, or strongly compact cardinals.

\bigskip

So we here see G\"odel expressing the view that axioms asserting the existence of inaccessible or Mahlo cardinals are ``intrinsically justified" in the sense of merely unfolding the content of the conception of set, whereas it was not clear at the time of his writing the 1966 footnote whether this was the case for axioms asserting the existence of weakly compact, measurable, or strongly compact cardinals, the only known justifications for these axioms at that time being ``extrinsic". An extrinsic justification for a candidate for a new axiom is based on something other than the idea that the axiom merely unfolds the content of the conception of set, such as for example the axiom having desirable or verifiable consequences.

\bigskip

G\"odel made the following remarks about this kind of justification:

\bigskip

\begin{quotation}
Secondly, however, even disregarding the intrinsic necessity of some new axiom, and even in case it has no intrinsic necessity at all, a probable decision about its truth is possible also in another way, namely, inductively by studying its ``success". Success here means fruitfulness in consequences, in particular in ``verifiable" consequences, i.e., consequences demonstrable without the new axiom, whose proofs with the help of the new axiom, however, are considerably simpler and easier to discover, and make it possible to contract into one proof many different proofs... There might exist axioms so abundant in their verifiable consequences, shedding so much light upon a whole field, and yielding such powerful methods for solving problems (and even solving them constructively, as far as that is possible) that, no matter whether or not they are intrinsically necessary, they would have to be accepted at least in the same sense as any well-established physical theory.' With the later work on Projective Determinacy in the 1980's, some suggested Projective Determinacy as a candidate for such an axiom.
\end{quotation}

\bigskip

These remarks were the origin of the distinction between ``intrinsic" and ``extrinsic" justifications which is still discussed by philosophers of set theory today.

\bigskip

In the rest of this work we shall be exploring various ways in which ``intrinsic justifications" for various reflection principles giving rise to large-cardinal axioms might be made out. Let us briefly review our philosophical ambitions.

\bigskip

We wish to argue for the following philosophical theses.

\bigskip

(1) In work of Schindler a challenge has been raised to the idea that any genuine intrinsic justification is available for any of the large-cardinal axioms, if we assume that we will only countenance sets and predicative proper classes. On this view, the usual attempts at intrinsic justification for the reflection principles that give rise to the usual large-cardinal axioms fail because higher-order quantification over the entire universe is not meaningful. Elaborating in somewhat more depth the ideas outlined in Tait's work on the ``bottom-up" conception, there may be an adequate response to this so that higher-order reflection principles may be intrinsically justified after all.

\bigskip

(2) There are natural strengthenings of Roberts-style reflection principles going at least up to the level of a supercompact cardinal.

\bigskip

(3) Intriguing parallels can be drawn between a certain natural family of weak reflection principles as outlined in previous work of the author and a different family similar to the set of principles presented by Roberts but with a very great jump upwards in consistency strength despite only minor modifications of the basic idea. This could lend some support to the view that Roberts-style reflection principles should indeed be seen as intrinsically justified, or alternatively perhaps one should say that there is some kind of ``bifurcation" of defensible views that could be put forward about the extent of intrinsic justification. On the two sides of this bifurcation, one of the available views would only allow for intrinsic justifications which do not go beyond the point of what is consistent with $V=L$, whereas the other available view would allow for intrinsic justifications for principles at the level of strength of those considered by Roberts.

\bigskip

The task of Section 2 shall be to argue for thesis 1. We shall explore what might be said about intrinsic justifications for small large cardinals such as inaccessible and Mahlo cardinals, which as noted above were indeed endorsed as intrinsically justified by G\"odel in his essay on the continuum problem. Thesis 1 will be defended in the course of that discussion. We shall also examine intrinsic justifications for reflection principles considered by Tait, Koellner, and the author, still at the level of what is consistent with $V=L$.

\bigskip

The task of Section 3 shall be to review contributions of Welch, Roberts and Marshall on reflection principles and intrinsic justifications for them, with a view towards defending thesis 2. The work of Marshall is not directly relevant to thesis 2, but it is of some value to explore it for the purpose of clarifying its relationship with the contribution of Roberts and stating some theorems. In addition we shall examine the relationship between reflection principles of similar strength to those considered by Roberts and weaker reflection principles discussed by the author, with a view towards defending thesis 3. Section 4 will consist of concluding remarks.

\bigskip

\section{Intrinsic justifications for small large cardinals}

In \cite{Tait2005a} Tait begins by quoting Cantor in \cite{Cantor1883} as saying, with reference to the totality $\Omega$ of transfinite ordinal numbers, `If the initial segment $\Sigma$ of $\Omega$ is a set, then it has a least strict upper bound $S(\Sigma)\in\Omega$.' Here we must understand `initial segment' so that initial segments are not necessarily proper. Thus we can see that it is obviously not consistent for every initial segment of $\Omega$ to be a set, since $\Omega$ itself cannot be a set. So some criterion must be given for when an initial segment is a set. Later Cantor attempted to solve this problem by invoking the notion of a `consistent multiplicity', but as Tait says `that is just naming the problem, not solving it.'  Tait proposes a modified version of Cantor's principle which Peter Koellner later in \cite{Koellner2009} called the Relativised Cantorian Principle, which says that for certain conditions $C$ for initial segments of $\Omega$ there corresponds an initial segment $\Omega_{C}$ such that `If the initial segment $\Sigma$ of $\Omega_{C}$ satisfies the condition $C$, then it has a least strict upper bound $S(\Sigma)\in\Omega_{C}$.' We shall spell this out in more detail with reference to a specific example below. With an appropriate choice of the condition $C$, which Tait calls an `existence condition', we can in this way determine an initial segment $\Omega_{C}$ of the totality of transfinite ordinal numbers $\Omega$. Tait writes `In this way, we are led to a hierarchy of more and more inclusive existence conditions, each of which can replace the condition ``\textit{is a set}" in Cantor's definition and so yields an initial segment of the transfinite numbers, but none of which yields ``all the numbers".' He calls this the `bottom-up' conception. The statement that we never obtain `all the numbers' by a procedure like this, that is, that we are indeed justified in supposing that $\Omega_{C}$ is a proper initial segment, might be thought to require further discussion; we shall consider this matter in more detail in due course.

\bigskip

We can articulate this notion of what is involved in this conception of set as follows. Let us begin with the example of motivating an inaccessible cardinal. The property of being inaccessible is equivalent to the property of being $\Pi^{0}_{n}$-indescribable for all positive integers $n$. Let us imagine that $V_{\alpha}$ is a level of the universe consisting of all the sets that we have built up so far, for example perhaps $\alpha$ is a worldly cardinal (a cardinal $\alpha$ such that $V_{\alpha}$ is a model of first-order ZFC). Now suppose we have some $X \subseteq V_{\alpha}$ and some formula $\phi$ in the second-order language of set theory with exactly one free second-order variable, and all quantified variables being first-order variables, with the property that the formula $\phi$, together with $X$ as the value assigned to its free second-order variable, witnesses the $\Pi^0_n$-describability of $V_{\alpha}$ for some positive integer $n$. The ordered pair $(\phi,X)$ constitutes a ``mark" for $V_{\alpha}$, which we regard as a warrant for concluding that $V_{\alpha}$ is not all of $V$. That is, the existence of such a mark is taken to be an ``existence condition" in the sense outlined by Tait as was described just above. We are therefore justified in adding one more level.

\bigskip

The idea that the existence of such a ``mark" is a warrant for concluding that the level $V_{\alpha}$ is not all of $V$ is based on the idea that $V$ is ``undefinable". The existence of a mark of the form described would count as a witness that the level built up so far is ``describable" and as such cannot be all of $V$. \footnote{I am thankful to William Tait for articulating this notion of a ``mark" in private email correspondence.} Now consider the totality of all sets that can be reached starting from the level $V_{0}$ by the process of adding one more level whenever we see a ``mark" for the level built up so far, and taking unions at limit ordinals. The totality of sets that can be built up in this way is a ``definable" totality, and so for that reason we should conclude that once again it cannot be all of $V$. Thus we are justified in positing a level of the universe $V_{\kappa}$ which is a closure point for the process of building up new levels in this way. These ideas serve as the motivation for higher-order reflection principles, using special instances of the Relativized Cantorian Principle.

\bigskip

It may be objected, along the line of thought outlined by Peter Koellner in \cite{Koellner2009}, that the Relativized Cantorian Principle proves too much, that assuming the existence of closure points for these kinds of processes begs the question. We certainly wouldn't regard a similar use of the Relativized Cantorian Principle to try to ``justify" a measurable cardinal to be a satisfactory argument without more detail being given. However in the case of the instances of the Relativized Cantorian Principle considered here we are motivating them by means of the idea that they are grounded in the ``undefinability" of $V$. We will later see examples of reflection principles that might be motivated in this way which appear to be \textit{prima facie} plausible but in fact turn out to be inconsistent. Thus we must acknowledge the reality that in seeking to justify reflection principles in this way we do need some appeal to ``extrinsic" evidence of consistency to ensure that we are not going astray. Nevertheless we can at least hold out the hope that further research into justifications of this form of reflection principles may yield insight into the content of the idea of the ``undefinability" of $V$ and enable us to flesh out a version of the iterative conception of sets on which the reflection principles in question can indeed be seen to be intrinsically justified in the sense of merely unfolding the content of the conception.

\bigskip

It may be instructive to consider a challenge posed by Schindler in \cite{Schindler1994} to the idea that any large cardinals are intrinsically justified. In that paper Schindler considers a class theory $BL_{1}$ which he claims to embody the weakest reflection principle schema that leads to the existence of large cardinals. He shows that it proves $\Delta^{1,NBG}_{1}$-comprehension but that this comprehension schema is false for the predicative classes. He concludes that we then face the dilemma between accepting the existence of non-predicative classes (for which he believes there is no good justification), or concluding that there is not sufficient justification for believing in the existence of large cardinals, at least in $V$ (which he finds to be an unsatisfying view). Can we say that the line of thought outlined above provides a good answer to Schindler's dilemma?

\bigskip

It may seem not, since it may seem that when I previously referred to the ``level of the universe" $V_{\alpha}$, for which I was contemplating whether or not to add one more level, I should in fact have been agnostic about whether it was a $V_{\alpha}$ or all of $V$, in order to avoid begging the question. If we are open to the possibility that it may be all of $V$ then it seems the suggested reflection principle based on the undefinability of $V$ is in fact false since we may have reached the point where we have all of $V$, and then $\Delta^{1,NBG}_{1}$-comprehension fails at that point, and this should count as a witness to ``definability", but we are not in fact justified in adding one more level since we now have all of $V$. On this account of the matter, I cannot claim to have answered Schindler's concerns without begging the question.

\bigskip

One possible way to answer this concern might be as follows. Yes, I do know that the collection of sets built up so far is indeed just a $V_{\alpha}$, since it was obtained using a certain circumscribed, definable set of procedures of formation, which as such cannot possibly exhaust all of $V$. So I am in a position to say that at this point I only have a $V_{\alpha}$ and impredicative comprehension holds. But I cannot use this fact as a ``warrant" for going further, since that would introduce an unacceptable self-referential quality into the definition of the phrase ``set of circumstances under which I say I have a warrant for going further", which would lead me to lose my justification for saying that proceeding in this way is not going to exhaust all of $V$. Although at any given stage we are in a position to say that we only have a set, the set of circumstances under which we say we have a warrant for going one level further must be more narrowly circumscribed than that. If this account of the matter is defensible, then it may be that we can make a claim to have answered Schindler's challenge and shown belief in some large cardinals to be warranted after all. This should only be taken as a tentative suggestion at this stage, but we may perhaps provisionally say that we have a possible response to Schindler's challenge to the idea that belief in large cardinals is justifiable.

\bigskip

The rest of the paper will examine various different types of reflection principles that one might seek to justify along these lines and the relations between them. This line of thought can be used to motivate the existence of levels of the universe satisfying higher-order reflection principles, assuming we accept a version of the iterative conception which is committed to the idea that ``the stages are undefinable". Of course, Boolos' influential exposition of the iterative conception in \cite{Boolos1971} did not view this idea as part of the iterative conception, and did not accept that even all of the Replacement axioms can be taken as part of the conception, but various other writers such as G\"{o}del in \cite{Godel1964} have taken something like such a requirement to be part of the content of the iterative conception. In \cite{Potter2004} Michael Potter gives an overview of the literature on this issue.

\bigskip

If we accept that higher-order reflection principles can be motivated in this way, the question then arises which reflection principles should be accepted. One naive idea might be as follows. Let $\mathcal{L}_{\in}^{\omega}$ be the $\omega$th-order language of set theory, and suppose we are considering a formula $\phi$ with free variables $x_1, x_2, \ldots x_{m_{1}}, X^{(2)}_{1}, X^{(2)}_{2}, \ldots X^{(2)}_{m_{2}}, X^{(3)}_{1}, X^{(3)}_{2}, \ldots \newline X^{(3)}_{m_{3}},  \ldots X^{(k)}_{1}, X^{(k)}_{2}, \ldots X^{(k)}_{m_{k}},$, where the variables denoted by lower-case letters are first-order variables, and for variables denoted by upper-case letters superscripts indicate the order of the variables, and an $m$th-order variable ranges over $\mathcal{P}^{m-1}(D)$ where $D$ is the domain of discourse. Suppose that the domain of discourse $D$ is equal to $V_{\kappa}$ for some ordinal $\kappa$. If $X^{(2)}\in\mathcal{P}(V_{\kappa})$, let $X^{(2),\alpha}:=X^{(2)}\cap V_{\alpha}$ for all ordinals $\alpha<\kappa$, and then inductively define $X^{(k),\alpha}=\{X^{(k-1),\alpha}\mid X^{(k-1)}\in X^{(k)}\}$ for all ordinals $\alpha<\kappa$ and all $X^{(k)}\in\mathcal{P}^{k-1}(V_{\kappa})$. Then we can posit that there should be a level $V_{\kappa}$ with the property that if the formula $\phi$ is true in $V_{\kappa}$ for a certain assignment of values to the variables, then there should be some ordinal $\alpha<\kappa$, such that all values assigned to free first-order variables are members of $V_{\alpha}$, and such that the relativisation of $\phi$ to $V_{\alpha}$ holds for an assignment of values to the variables such that if a variable has $X^{(j)}$ assigned to it under the first assignment then it will have $X^{(j),\alpha}$ as defined above assigned to it under the second assignment. By the relativisation of $\phi$ to $V_{\alpha}$ we mean the result of relativising the first-order quantifiers to $V_{\alpha}$, the second-order quantifiers to $V_{\alpha+1}$, and so on. This is a naive form of reflection which one might naturally posit.

\bigskip

If one makes the restriction that the free variables should be no higher than second-order then this is consistent (relative to the consistency of a totally indescribable cardinal), but if one allows parameters of third order or higher then this is inconsistent, as was first observed by Reinhardt in \cite{Reinhardt1974}, but can also be seen as follows, using an argument presented by Tait in \cite{Tait2005a}, although a proof can also be found on p. 274 of Drake in \cite{Drake1974}. Let $\phi$ be a formula with exactly one free third-order variable $X^{(3)}$ asserting that every element of $X^{(3)}$ is a bounded set of ordinals. Then if we assign to the variable $X^{(3)}$ the set $\{\{\alpha\mid\alpha<\beta\}\mid\beta<\kappa\}$, then $\phi(X^{(3)})$ will be true in $V_{\kappa}$ but we will not have $\phi^{\alpha}(X^{(3),\alpha})$ for any $\alpha<\kappa$, where $\phi^{\alpha}$ denotes the relativisation of $\phi$ to $V_{\alpha}$. Thus the reflection principle described is inconsistent if we allow parameters of third order or higher. So if one wants to accept the reflection principle in the case of parameters no higher than second order, one must identify principled reasons why one should not also accept it for parameters of third order or higher.

\bigskip

Let us consider what Tait writes in \cite{Tait2005a} about this question.

\bigskip

`One plausible way to think about the difference between reflecting $\phi(A)$ when $A$ is second-order and when it is of higher-order is that, in the former case, reflection is asserting that, if $\phi(A)$ holds in the structure $\langle R(\kappa),\in,A \rangle$, then it holds in the substructure $\langle R(\beta),\in, A^{\beta} \rangle$ for some $\beta<\kappa \ldots$ But, when $A$ is higher-order, say of third-order, this is no longer so. Now we are considering the structure $\langle R(\kappa), R(\kappa+1), \newline \in, A \rangle$ and $\langle R(\beta), R(\beta+1), \in, A^{\beta} \rangle$. But, the latter is not a substructure of the former, that is the `inclusion map' of the latter structure into the former is no longer single-valued: for subclasses $X$ and $Y$ of $R(\kappa), X\neq Y$ does not imply $X^{\beta}\neq Y^{\beta}$. Likewise for $X\in R(\kappa+1), X\notin A$ does not imply $X^{\beta}\notin A^{\beta}$. For this reason, the formulas that we can expect to be preserved in passing from the former structure to the latter must be suitably restricted and, in particular, should not contain the relation $\notin$ between second- and third-order objects or the relation $\neq$ between second-order objects.'

\bigskip

He then uses these ideas to motivate the following family of reflection principles.

\begin{Definition} A formula in the $\omega$th-order language of set theory is positive iff it is built up by means of the operations $\vee$, $\wedge$, $\forall$, $\exists$ from atoms of the form $x=y$, $x \neq y$, $x \in y$, $x \notin y$, $x \in Y^{(2)}$, $x \notin Y^{(2)}$ and $X^{(m)} = X'^{(m)}$ and
$X^{(m)} \in Y^{(m+1)}$, where $m \geq 2$. \end{Definition}

\begin{Definition} Suppose that $A$ is a variable of order at most 2. Define $A_{n\times}=\{\langle n,x\rangle\mid x\in A\}$ and $A_{/n}=\{x\mid\langle n,x\rangle\in A\}$. In the case where $B$ is a variable of order greater than 2, define, using induction on the order of $B$, $B_{n\times}:=\{A_{n\times}\mid A\in B\}$ and $B_{/n}:=\{A_{/n}\mid A\in B\}$. Given an $n$-tuple of variables of order $m$, $(X_{0}^{(m)}, X_{1}^{(m)}, \ldots X_{n-1}^{(m)})$, call the operator $S$ which sends this $n$-tuple to $S(X_{0}^{(m)}, X_{1}^{(m)}, \ldots X_{n-1}^{(m)}):=\bigcup_{i=0}^{n-1} (X_{i}^{(m)})_{i\times}$ a contracting operator, and call an operator of the form $X \mapsto X_{/n}$ a de-contracting operator. Call a formula $\phi$ in the $\omega$th-order language of set theory positive in the extended sense if it is equivalent to a formula, which has no quantified variables of order higher than first-order, and is positive, except that it may include contracting and de-contracting operators. \end{Definition}

\begin{Definition} For $0<n<\omega$, $\Gamma^{(2)}_{n}$ is the class of formulas

$$ \forall X_{1}^{(2)} \exists Y_{1}^{(k_{1})} \cdots \forall X_{n}^{(2)} \exists Y_{n}^{(k_{n})} \phi(X_{1}^{(2)}, Y_{1}^{(k_{1})}, \ldots, X_{n}^{(2)}, Y_{n}^{(k_{n})}, A^{(l_{1})}, \ldots A^{(l_{n'})}) $$

where $\phi$ is positive in the extended sense and $k_{1}, \ldots k_{n}, l_{1}, \ldots l_{n'}$ are natural numbers which are greater than or equal to 2. \end{Definition}

\begin{Definition} We say that $V_{\kappa}$ satisfies $\Gamma^{(2)}_{n}$-reflection if, for all $\phi\in\Gamma^{(2)}_{n}$, if \newline $V_{\kappa}\models\phi(A^{(m_{1})},A^{(m_{2})},\ldots A^{(m_{p})})$ then $V_{\delta}\models\phi(A^{(m_{1}),\delta},A^{(m_{2}),\delta},\ldots A^{(m_{p}),\delta})$ for some $\delta<\kappa$, where we relativize the higher-order parameters in the manner described previously. \end{Definition}

Peter Koellner established in \cite{Koellner2009} that these reflection principles are consistent relative to an $\omega$-Erd\H{o}s cardinal. In \cite{Tait2005a} Tait proposes to define $\Gamma^{(m)}_{n}$ in the same way as the class of formulas $\Gamma^{(2)}_{n}$, except that universal quantifiers of order $\leq m$ are permitted. Koellner shows in \cite{Koellner2009} that this form of reflection is inconsistent when $m>2$.

\bigskip

Now I shall describe a family of reflection principles which are consistent relative to an $\omega$-Erd\H{o}s cardinal and which yield all of the reflection principles of Tait considered so far not known to be inconsistent. I shall describe a naive version of it first which is inconsistent and then show how to modify it so as to make it consistent relative to an $\omega$-Erd\H{o}s cardinal.

\bigskip

Consider a formula $\phi$ in the $\omega$th-order language of set theory with free variables $x_{1}, x_{2}, \ldots x_{m_{1}}, X^{(2)}_{1}, X^{(2)}_{2}, \ldots X^{(2)}_{m_{2}}, \ldots X^{(k)}_{1}, X^{(k)}_{2}, \ldots X^{(k)}_{m_{k}}$ and suppose that $\phi$ holds in $V_{\kappa}$ for a certain assignment of values to the free variables. The first idea we try is simply to avoid re-interpreting the higher-order variables when we reflect downwards. That is, we posit that there should be an $\alpha<\kappa$ with $x_{1}, x_{2}, \ldots x_{m_{1}} \in V_{\alpha}$, such that $\phi$ holds with all of the free variables interpreted the same way as before, first-order variables ranging over $V_{\alpha}$ and $k$th-order variables ranging over $\mathcal{P}^{k-1}(V_{\kappa})$ for $k\geq 2$. This is a natural form of reflection to consider as an alterative to the first family of reflection principles we considered that was shown to be inconsistent. However, this form of reflection is also unfortunately inconsistent. To see this, consider the sentence $\phi$ given by $\forall X \forall Y (X \neq Y \implies \exists x ((x \in X \wedge x \notin Y) \vee (x \notin X \wedge x \in Y)))$, where $X$ and $Y$ are second-order variables. This formula clearly gives a counter-example to our reflection principle for any level $V_{\kappa}$.

\bigskip

We can diagnose the root of this problem as arising from the fact that we introduced an existential quantifier for a first-order variable within the scope of a higher-order quantifier. This means that when we ``Skolemize" the formula the existentially quantified first-order variable becomes a Skolem function of a tuple of variables including higher-order variables, and now the fact that the number of possible values for a higher-order variable is large may mean that we cannot have closure under the Skolem functions for any reflecting structure $(V_{\alpha}, V_{\kappa+1}, V_{\kappa+2}, \ldots)$ for any $\alpha<\kappa$. If we modify the reflection principle so as only to apply to formulas in which unbounded existential quantifiers for first-order variables do not appear within the scope of a higher-order quantifier (but existential quantifiers bound by a first-order variable are allowed), then the resulting reflection principle becomes consistent relative to an $\omega$-Erd\H{o}s cardinal. In fact, as shown in \cite{McCallum2020} (see in particular Lemma 1.7), a level $V_{\kappa}$ satisfies this form of reflection if and only if $\kappa$ is $\omega$-reflective in the sense defined in \cite{McCallum2013}, and when this is so then $V_{\kappa}$ satisfies $\Gamma^{(2)}_{n}$-reflection for all $n \in\omega$.

\begin{Definition} \label{reflection} The theory in the $\omega$th-order language of set theory outlined in the paragraph above shall be denoted by $T$. \end{Definition}

It is natural to posit this as the correct modification of the inconsistent reflection principle considered earlier.

\begin{Definition} If $n\in\omega$ and $n>0$, a cardinal $\kappa$ is said to be $n$-ineffable if for every function $f:[\kappa]^{n+1}\rightarrow 2$ there is a stationary subset $S\subseteq \kappa$ which is homogeneous for $f$, that is, $f$ is constant on $[S]^{n+1}$. A cardinal $\kappa$ is said to be totally ineffable if it is $n$-ineffable for every $n\in\omega\setminus\{0\}$.

\end{Definition}

An $\omega$-reflective cardinal is totally ineffable, as follows from results proved in \cite{McCallum2013}, and every totally ineffable cardinal is a stationary limit of totally indescribable cardinals, so we obtain a justification for the form of reflection considered earlier in the case of parameters no higher than second order.

\bigskip

This completes our discussion of intrinsic justifications for small large cardinals. It may be well to make a few remarks clarifying why these suggested justifications are ``intrinsic". We are basing this all on features of the concept of set. We are not appealing to any results from inner model theory or any ``extrinsic" justifications of the type considered in Section 3 of \cite{Koellner2006}. Thus these count as ``intrinsic" justifications in the sense outlined in the Introduction. All cardinals considered so far are indeed small in the sense of being consistent relative to an $\omega$-Erd\H{o}s cardinal. The time has now come to consider the reflection principles of \cite{Roberts2017} and \cite{Marshall1989}.

\section{The Reflection Principles of Welch, Roberts and Marshall}

Let us now briefly describe the reflection principle introduced by Roberts in \cite{Roberts2017}. To explain the reflection principle which Roberts formulates in this paper, let us begin by explaining the reflection principle that he calls $\textsf{R}_{2}$. This is an axiom schema in the second-order language of set theory. For each formula $\phi(x_{1}, x_{2}, \ldots x_{m}, X_{1}, X_{2}, \ldots X_{n})$ in the second-order language of set theory, there is an axiom asserting that if $\phi$ holds, then there exists an ordinal $\alpha$ such that $x_{1}, x_{2}, \ldots x_{m} \in V_{\alpha}$, and a ``set-sized" family of classes which contains the classes $X_{1}, X_{2}, \ldots X_{n}$, which is itself coded for by a single class, and which is standard for $V_{\alpha}$ in the sense that every subset $X\subseteq V_{\alpha}$ is such that some class in the family has intersection with $V_{\alpha}$ equal to $X$, such that the formula $\phi$ still holds when the first-order variables are relativised to $V_{\alpha}$ and the second-order variables are relativised to the set-sized family of classes. This completes the description of the axiom schema $\textsf{R}_{2}$. Then Roberts extends the axiom schema as follows. He extends the underlying language so as to include a satisfaction predicate for the second-order language of set theory, and then he extends the axiom schema so as to also include an axiom of the kind described for every formula in this extended language, calling this new axiom schema $\textsf{R}_{S}$. Then he denotes by $\textsf{ZFC2}_{S}$ the result of extending $\textsf{ZFC2}$ -- being the same as ZFC except for having Separation and Replacement as single second-order axioms and also having an axiom schema of class comprehension for every formula in the second-order language of set theory -- by adding the usual Tarskian axioms for the satisfaction predicate and extending the class comprehension axiom schema to include axioms involving formulas in the extended language. Then he proceeds to investigate the theory $\textsf{ZFC2}_{S}+\textsf{R}_{S}$. He investigates the version of the reflection principle for this particular language in order to obtain the degree of large-cardinal strength that he wants, but clearly the basic underlying idea can be applied to a broad class of languages. In particular, regarding $\textsf{ZFC2}+\textsf{R}_{2}$, this theory is consistent relative to a $\Pi^{1}_{\infty}$-indescribable cardinal, thus being relatively weak, but $\textsf{ZFC2}_{S}+\textsf{R}_{S}$ is considerably stronger as we shall clarify presently, and later we shall go on to consider other variants applying the same basic idea to more expressive languages thereby yielding theories that are stronger still.

\bigskip

This completes the description of the reflection principle which Roberts considers. He shows that the theory $\textsf{ZFC2}_{S}+\textsf{R}_{S}$ proves the existence of a proper class of 1-extendible cardinals and is consistent relative to a 2-extendible cardinal. The reflection principle suggested by Welch in \cite{Welch2017} is similar, and makes use of an elementary embedding rather than a ``reflecting structure".

\bigskip

We can note the similarity between Roberts' reflection principle and the reflection principle described in the previous section. Specifically, both reflection principles avoid any re-interpretation of parameters, and involve ``reflecting structures" in which the domain of the first-order variables is reflected to a lower rank but the higher-order variables are still interpreted by entities of the same rank as before.

\bigskip

To further illuminate the discussion, it will be useful to analyze further what kinds of comparisons can be made between the reflection principle discussed in the previous section and the reflection principle of Roberts. But, before doing so, we shall introduce the reflection principles discussed by Marshall in \cite{Marshall1989} and discuss how they are related to the reflection principle of Roberts (and natural generalizations of it).

\bigskip

So let us now describe Marshall's family of theories $B_n$ indexed by the positive integers $n$, following her exposition in \cite{Marshall1989}. The first theory in the sequence, $B_1$, was formulated by Bernays and is presented in \cite{Chuaqui1978}. We shall first describe the theory $B_1$ and then proceed to describe Marshall's generalizations $B_n$ for integers $n>1$.

\bigskip

The theory $B_1$ is a theory in the first-order language of set theory. The objects of the domain of discourse are called classes and a class is said to be a set if and only if it is a member of some class, otherwise it is said to be a proper class. Now we shall describe the axioms of the theory $B_1$.

\bigskip

(A1) Extensionality. $\forall x (x \in A \equiv x \in B) \implies A=B$.

\bigskip

(A2) Class specification. Suppose that $\phi$ is a formula and that $A$ is not free in $\phi$. Then

$$ \exists A \forall x (x \in A \equiv \phi(x) \wedge x \hspace{1 mm} \mathrm{is} \hspace{1 mm} \mathrm{a} \hspace{1 mm} \mathrm{set}). $$

\bigskip

(A3) Subsets. $A \hspace{1 mm} \mathrm{is} \hspace{1 mm} \mathrm{a} \hspace{1 mm} \mathrm{set} \wedge B\subseteq A \implies B \hspace{1 mm} \mathrm{is} \hspace{1 mm} \mathrm{a} \hspace{1 mm} \mathrm{set}.$

\bigskip

(A4) Reflection principle. Suppose that $\phi(x)$ is a formula. Then

$$ \phi(A) \implies \exists u (u \hspace{1 mm} \mathrm{is} \hspace{1 mm} \mathrm{a} \hspace{1 mm} \mathrm{transitive} \hspace{1 mm} \mathrm{set} \wedge \phi^{\mathcal{P}u}(A \cap u)). $$

\bigskip

(A5) Foundation. $\exists x (x \in A) \implies \exists x (x \in A \wedge \forall z (z \notin x \vee z \notin A))$.

\bigskip

(A6) Choice for sets. $\forall x((x \hspace{1 mm} \mathrm{is} \hspace{1 mm} \mathrm{a} \hspace{1 mm} \mathrm{set} \wedge \forall y \forall y' ((y \in x \wedge y' \in x \wedge y \neq y') \implies y\cap y'=\emptyset)) \implies \exists x' (\forall y (y \in x \implies \exists ! z (z \in y \wedge z \in x'))))$.

\bigskip

This completes the description of the axioms of the theory $B_1$.

\begin{theorem} The natural models of $B_1$ are those sets of the form $V_{\alpha+1}$ for some ordinal $\alpha$ such that $(V_{\alpha+1},\in)$ is a model for the theory $B_1$. An ordinal $\alpha$ is such that $V_{\alpha+1}$ is such a natural model if and only if $\alpha$ is a $\Pi^1_\infty$-indescribable cardinal. \end{theorem}

\begin{proof}

It is easily seen that every $V_{\alpha+1}$ is a model of (A1)-(A3) and (A5)-(A6) so therefore this is simply a question of noting that satisfaction of (A4) by $V_{\alpha+1}$ can be seen to be equivalent to $\alpha$ being $\Pi^1_\infty$-indescribable. This equivalence can be seen, when we consider that range of the variable $u$ is over all transitive sets, but the truth-value of the assertion is not changed if we restrict to transitive sets of the form $V_{\beta}$.

\end{proof}

From this theorem stated immediately above, we see that the existence of a $\Pi^1_\infty$-indescribable cardinal is a sufficient assumption on which to prove that $B_1$ is consistent. The theory $B_1$ in turn implies the relativization of all the axioms of $\textsf{ZFC}$ to the class of sets (as shown in \cite{Marshall1989}) and also (the relativization to the class of sets of) the existence of a proper class of $\Pi^1_n$-indescribable cardinals for each positive integer $n$. Let us now see how Marshall generalized this to a sequence of stronger theories. We shall first describe the theory $B_2$ and then proceed to describe the theories $B_n$ for $n>2$.

\bigskip

The theory $B_2$ is a theory in the first-order language of set theory. The objects of the domain of discourse are called 2-classes, an object is said to be a class if and only if it is a member of some 2-class, and an object is said to be a set if and only if it is a member of some class. The class of all sets (whose existence and uniqueness will follow from the axioms soon to be presented) is denoted by $V_0$, and the 2-class of all classes is denoted by $V_1$.

\bigskip

The first axiom for the theory is the axiom of 2-class specification.

\bigskip

(A1) Suppose that $\phi$ is a formula and $A$ is not free in $\phi$. Then

$$ \exists A \forall x (x \in A \equiv \phi(x) \wedge x \hspace{1 mm} \mathrm{is} \hspace{1 mm} \mathrm{a} \hspace{1mm}  \mathrm{class}). $$

Next we have extensionality for 2-classes.

\bigskip

(A2) $\forall x (x \in A \equiv x \in B) \implies A=B$.

\bigskip

We write $\mathrm{STC}(v)$ to abbreviate a formula saying that $v$ is transitive, and every sub-2-class $x$ of any element $y$ of $v$ is also an element of $v$.

\bigskip

If $A$ and $u$ are 2-classes, $A^u$ denotes $A \cap u$ if $A$ is a class, and denotes $\{x^u\mid x \in A \cap u\}$ if $A$ is not a class. $\phi^u$ denotes the relativization of the formula $\phi$ to $u$.

\bigskip

If $x$ and $y$ are two classes, then the ordered pair $[x,y]$ of the classes $x$ and $y$ is defined to be $x \times \{0\} \cup y \times \{1\}$. The notation $A(u)$ is an abbreviation for $u \cap V_0 \in V_0 \wedge \mathrm{STC}(u \cap V_0) \wedge \forall x \forall y (x, y \in u \implies [x,y] \in u)$. Then our reflection principle is given by

\bigskip

(A3) If $\phi(x)$ is a formula, then

$$ \phi(A) \implies \exists u (A(u) \wedge \phi^{\mathcal{PP}(u \cap V_0)}(A^u)). $$

Finally we take as axioms the axiom of foundation and the axiom of choice for sets as in the theory $B_1$, noting that the interpretation of the abbreviation ``$x$ is a set" is now different. This completes the description of the theory $B_2$. Given any theorem $\phi$ of $B_1$, the relativization of $\phi$ to $V_1$ is a theorem of $B_2$, as shown in \cite{Marshall1989}. Marshall also shows in \cite{Marshall1989} that $B_2$ proves that the class of 1-extendible cardinals is a stationary subclass of the class of ordinals, and that if $\kappa$ is 2-extendible or $\beth_{\kappa+1}$-supercompact, then $V_{\kappa+2}$ is a natural model of $B_2$. Before proceeding to give the description of the theories $B_n$ with $n>2$, let us pause to discuss the relationship between $B_2$ and a natural generalization of Roberts' theory $\textsf{ZFC2}_{S}+\textsf{R}_{S}$, and also to characterize the natural models of $B_2$ and related theories, thereby answering a question asked by Marshall in \cite{Marshall1989}.

\bigskip

There is an obvious extension of Roberts' theory $\textsf{ZFC2}_{S}+\textsf{R}_{S}$ to a theory in the third-order language of set theory $\textsf{ZFC3}+\textsf{R}_{3}$. And there is an obvious interpretation of the third-order language of set theory in the first-order language of set theory where the domain of discourse of the variables in the latter language is understood to be the 2-classes.

\begin{theorem} The image of the theory $\textsf{ZFC3}+\textsf{R}_{3}$ under this interpretation is equivalent to the theory $B_2$. \end{theorem}

\begin{proof} In the theory $\textsf{ZFC3}+\textsf{R}_{3}$ a form of reflection is posited in which the reflecting structure has first-order variables ranging over some $V_{\alpha}$ and higher-order variables ranging over some ``set-sized" family of classes and 2-classes which is standard for $V_{\alpha}$. The difference in the case of $B_2$ is that firstly the first-order variables range over the extension of some $u$ which may or may not be a $V_{\alpha}$, and the higher-order variables are reflected downward as well, so as to range over $\mathcal{P}(u)$ and $\mathcal{P}(\mathcal{P}(u))$ respectively. The latter issue, that the higher-order variables are reflected downward, will clearly make no difference if the ``set-sized" family of classes and 2-classes which is posited as the reflecting structure by $\textsf{ZFC3}+\textsf{R}_{3}$ can be assumed without loss of generality to satisfy the appropriate forms of extensionality, and clearly this is the case for $\textsf{ZFC3}+\textsf{R}_{3}$. The fact that possibly the first-order variables may be reflected to a $u$ rather than a $V_{\alpha}$ makes no difference since clearly the theory $B_2$ entails that it can be assumed without loss of generality to be a $V_{\alpha}$. Thus from these considerations we see that the image of $\textsf{ZFC3}+\textsf{R}_{3}$ under the appropriate interpretation into the first-order language of set theory is indeed equivalent to $B_2$ as claimed. \end{proof}

\bigskip

Similar relationships can be found between other natural extensions of Roberts' theory and the theories $B_n$ for integers $n>2$ to be defined later. Let us characterize the natural models of $B_2$. Roberts shows in \cite{Roberts2017} that if $V_{\kappa+2}$ is a natural model of $B_2$ then $\kappa$ is a limit of a 1-extendible chain, that is, a subset $S \subseteq \kappa$ cofinal in $\kappa$ such that if $\alpha \in S, \alpha<\beta$, and $\beta \in S \cup \{\kappa\}$, then $\alpha$ is 1-extendible to $\beta$, that is, there is an elementary embedding $j: V_{\alpha+1} \prec V_{\beta+1}$ with critical point $\alpha$ and $j(\alpha)=\beta$. Given any positive integer $n$ there is a clear corresponding notion of a $\Pi^2_n$-extendible chain, and Roberts' result easily generalizes to show that $V_{\kappa+2}$ is a natural model of $B_2$ only if $\kappa$ is a limit of a $\Pi^2_n$-extendible chain for all positive integers $n$. But the converse of this is also easy to prove; one easily sees that a $\kappa$ which is a limit of a $\Pi^2_n$-extendible chain for all positive integers $n$ will be such that $V_{\kappa+2}$ is a natural model of $B_2$. Thus we have obtained a characterization of the natural models of $B_2$.

\bigskip

There is a similar characterization of the natural models of $\textsf{ZFC2}_{S}+R_{S}$ in which $V_{\alpha+2}$ is replaced with $L_{1}(V_{\alpha+1})$ (by which we mean the definable powerset of $V_{\alpha+1}$) and one uses the hierarchy of complexity classes of formulas relativized to $L_{1}(V_{\alpha+1})$. Let us now return to the task of defining the theories $B_n$ for integers $n>2$.

\bigskip

The theory $B_n$ for an integer $n>2$ is a theory in the first-order language of set theory where the domain of discourse is taken to be the $n$-classes, where classes are collections of sets, 2-classes are collections of classes, 3-classes are collections of 2-classes, and so on. The axioms of extensionality and class specification are as before, and there is an obvious definition of $V_0$, the collection of sets, $V_1$, the collection of classes, and $V_k$, the collection of $k$-classes for an integer $k$ such that $2 \leq k \leq n$. For arbitrary $n$-classes $A$ and $u$, we define $A^u=A \cap u$ if $A \in V_1$ and then by $\in$-induction we define $A^u=\{x^u\mid x \in A \cap u\}$ for $A \notin V_1$. Then our reflection principle is: For each formula $\phi(x)$, we have

$$ \phi(A) \implies \exists u (A(u) \wedge \phi^{\mathcal{P}^n(V_0 \cap u)}(A^u)). $$

Then we have the axiom of foundation and the axiom of choice for sets as before. This completes the definition of the theory $B_n$, which is equivalent to natural extensions of Roberts' theory $\textsf{ZFC2}_{S}+\textsf{R}_{S}$ in the same way as before.

\bigskip

Suppose that we consider an extension of Roberts' theory to a theory in the $\omega$th-order language of set theory which we might label $\textsf{ZFC}\omega+\textsf{R}_{\omega}$, and suppose that we drop the requirement that the reflecting structure must be ``set-sized", but keep the requirement that it must be ``standard". Saying that the reflecting structure has to be ``set-sized" amounts to saying that the domains of variables of all types are strictly smaller than $\beth_{\kappa}$ where $\kappa$ is the cardinal being reflected, and requiring it to be ``standard" means that the image of the domain of each type of higher-order variable in the reflecting structure (a subset of $V_{\kappa+n-1}$ where $n$ is the type of the variable) under the collapsing map to $V_{\alpha+n-1}$ should be full in the sense of being all of $V_{\alpha+n-1}$. This modification does not change the class of natural models for the theory, and the natural models of this theory are of the form $V_{\kappa+\omega}$ where $\kappa$ is a limit of an $n$-extendible chain for each positive integer $n$. We could also easily formulate an essentially equivalent theory in the first-order language of set theory along the lines outlined by Marshall, where rather than having a ``reflecting structure" which fails to re-interpret higher-order parameters downward, we have a reflecting structure which serves to ``guide" the reflection of the higher-order parameters. Here we are alluding to Marshall's notation $A^u$, explained in her paper, which explains how the structure $u$, taking the place of what for Sam Roberts is the reflecting structure, serves to determine the parameter $A^u$ to which the higher-order parameter $A$ is reflected. Once again such a variant of Marshall's theories would be essentially equivalent to $\textsf{ZFC}\omega+\textsf{R}_{\omega}$.

\bigskip

On the other hand, I might modify the reflection principle in the $\omega$th-order language of set theory described in Section 2, so that the reflecting structure is allowed to be any structure of the form $\langle V_{\alpha}, X_1, X_2, \ldots \rangle$ where $\alpha<\kappa$ and $X_i \subseteq V_{\kappa+i}$ for $i\in\omega$ and $i>0$, and each $X_i$ is standard for $V_{\alpha}$ in the sense that the image under the transitive collapsing map contains all of $V_{\alpha+i}$. That is, the permitted reflecting structures are exactly the same as those allowed in the previously described modification of Roberts' theory, rather than those with a full domain for the higher-order parameters as was the case in Section 2. The resulting theory, which we shall denote by $T'$ (see statement of theorem below), is in fact equivalent to the theory $T$ originally described in Section 2 in Definition \ref{reflection}, and the natural models of $T$ still correspond in the obvious way to the $\omega$-reflective cardinals. I must now provide a citation which would enable one to re-construct an argument for that point.

\begin{theorem} Let $T$ be as in Definition \ref{reflection}, and let us denote by $T'$ the theory described above. The theories $T$ and $T'$ are equivalent and the natural models of $T$ correspond in the obvious way to the $\omega$-reflective cardinals. \end{theorem}

\begin{proof} A slight modification of the argument given for Lemma 1.7 in \cite{McCallum2020} is sufficient to establish the claim. \end{proof}

\noindent The only difference between the theory $T$ described above and the modification of Roberts' theory $\textsf{ZFC}\omega+\textsf{R}_{\omega}$ is that the theory $T$ does not allow reflection of formulas which have unbounded existential quantification of first-order variables within the scope of a quantifier for a higher-order variable. However, dropping that requirement on the formulas only leads to a consistent reflection principle if we allow reflecting structures which are ``standard" rather than ``full". Thus we have clarified the relationship between the theory $T$ (and variants thereof) and the theory $\textsf{ZFC}\omega+\textsf{R}_{\omega}$, and have also perhaps provided some motivation for the view that if the former is intrinsically justified then so is the latter.

\bigskip

This indicates the kind of comparison that can be made between the reflection principles that give rise to $\omega$-reflective cardinals, and the much stronger reflection principles considered by Roberts and Marshall which give rise to $n$-extendible cardinals for every positive integer $n$. That is, given an appropriate choice of formulation of each reflection principle in question, they are identical except with regard to whether or not unbounded existential quantifiers for first-order variables within the scope of a higher-order quantifier are to be permitted in the formulas being reflected.

\bigskip

We shall discuss in some detail below how to construct a reflection principle stronger than that of Roberts, which stands in the same relation to the reflection principle which is shown in \cite{McCallum2020} to be equivalent to a remarkable cardinal, as does the just recently mentioned stronger reflection principle giving rise to $n$-extendible cardinals for all $n$, to the $\omega$-reflective cardinals, and it will turn out that the large-cardinal property arising from this reflection principle is in fact equivalent to supercompactness. Let us briefly recall a slight variant on the definition of the notion of extremely reflective cardinal given in \cite{McCallum2020}, equivalent to the definition actually given there.

\begin{Definition} A cardinal $\kappa$ is said to be extremely reflective if, for each ordinal $\lambda>\kappa$, considering structures of the form $(V_{\kappa}, V_{\lambda}\setminus V_{\kappa})$ and formulas $\phi$ in a two-sorted language holding in such a structure, with the requirement that only bounded existential quantifiers for variables of the first sort should appear within the scope of quantifiers for variables of the second sort, each such formula reflects down to some $\beta<\kappa$ in the sense that any first-order parameters appear in $V_{\beta}$, the reflecting structure is of the form $(V_{\beta}, V_{\lambda}\setminus V_{\kappa})$, and all parameters corresponding to variables of the second sort remain as before, not ``reflected downward". \end{Definition}

Clearly one could naturally modify the reflection principle described above by modifying the range of allowed reflecting structures to include those which are ``standard", in the sense that the image of the range of the second sort of variables in the reflecting structure under the collapsing map (which need not be injective) includes all of $V_{\gamma} \setminus V_{\beta}$ for some $\gamma$ such that $\beta<\gamma<\kappa$, and ``set-sized" in the sense that the cardinality of the domain of the variables of second sort in the reflecting structure is less than $\beth_{\kappa}$, but not ``full" in the sense that the variables of second sort range over all of $V_{\lambda} \setminus V_{\kappa}$. We will also drop the constraint on quantifiers in the formulas that are reflected. Once again we have a weak version of the reflection principle, given above, and proved in \cite{McCallum2020} to be equivalent to a remarkable cardinal (see in particular Theorem 2.3 of the cited work), and a stronger version, outlined just now, which we shall presently prove is equivalent to a supercompact cardinal.

\bigskip

In this way we can see two different possible notions of intrinsic justification, one extending all the way up to a remarkable cardinal, and one extending all the way up to a supercompact cardinal.

\bigskip

\begin{theorem} $V_{\kappa}$ satisfies the stronger form of reflection outlined above if and only if $\kappa$ is supercompact. \end{theorem}

\begin{proof} The point can be argued by means of Magidor's characterisation of supercompactness in terms of elementary embeddings between ranks, a proof of which can be found in \cite{Kanamori2003}. It is basically just a question of confirming that the satisfaction by $V_{\kappa}$ of the reflection principle in question is equivalent to the existence of an elementary embedding of the desired kind.

\bigskip

It may be instructive to compare with the proof given in \cite{McCallum2020} that a cardinal is extremely reflective if and only if it is remarkable. In \cite{McCallum2020} a cardinal was said to be extremely reflective if it satisfied the weaker version of the reflection principle described above, and this was shown to entail remarkability by deriving the existence of an elementary embedding $j:(V_{\alpha})^{V} \rightarrow (V_{\beta})^{V}$ for appropriately chosen $\alpha, \beta$, where $j \in V[G]$ for some generic extension $V[G]$ of $V$. (This was done by means of a result of Schindler showing that the existence of such a ``virtual" elementary embedding is equivalent to the existence of a winning strategy in a particular game.)

\bigskip

The argument for why the stronger reflection principle entails supercompactness is simpler, since we do not need to deal with ``virtual" elementary embeddings or game-theoretic characterizations thereof. Rather, we are seeking an elementary embedding in the ground model $V$, and it is possible to recover it from the inverse of the collapsing map applied to the range of the variables of the second sort in the reflecting structure, in the case where the formula reflected includes extensionality as one of its conjuncts, as we may assume without loss of generality. (In the general case it was remarked that the collapsing map need not be injective, but is clearly is so in the case where the formula being reflected includes extensionality as a conjunct.) It may be complained that this construction of an elementary embedding will only achieve elementarity ``one formula at a time", but $\lambda$ was arbitrary so we may obtain elementarity for a truth predicate for $V_{\lambda}$ by means of a formula relativized to $V_{\lambda+1}$ together with a ``truth set" as parameter. Thus we quite easily obtain that the stronger reflection principle for $\kappa$ implies that $\kappa$ is supercompact as characterized in terms of elementary embeddings by Magidor.

\bigskip

Thus in the case of the stronger reflection principle the argument proceeds much more simply than in the case of the weaker reflection principle, discussed in \cite{McCallum2020}, where it is shown in detail how to construct a ``virtual" elementary embedding in that case. It may still require some discussion why the converse holds, that is why supercompactness entails the stronger version of the reflection principle. This is basically a question of showing that, if we assume that the elementary embeddings called for in Magidor's characterization of supercompactness are available, then they can be chosen without loss of generality to have arbitrary parameters from $V_{\lambda}$ in their range. This can be established by means of a Skolem hull argument.

\end{proof}

This form of reflection is quite closely related to, but not quite as strong as, the form of reflection considered in Marshall's theory $B_0(V_0)$ of the paper \cite{Marshall1989}. This theory is structured similarly to the theories $B_n$, but it has a universe $V_0$ embedded in the larger universe $V$ which is assumed to be a model of ZFC, and the reflection principle reflects formulas holding in $V$ in such a way that $V_0$ is reflected to $V_{\kappa}$ for some $\kappa \in V_0$ and $V$ is reflected to $V_0$. The reflection of the parameters is relative to $V_{\kappa}$ and similar to as in the theories $B_n$. This is clearly a stronger form of reflection than the one already considered. With the reflection principle described above, $V_{\kappa}$ and $V_{\lambda}$ were being reflected to $V_{\alpha}$ and $V_{\beta}$ respectively with $\alpha, \beta < \kappa$, and the form of reflection in $B_0(V_0)$ is similar except that $V_{\lambda}$ is replaced by $V$ and $V_{\beta}$ by $V_{\kappa}$. Every instance of the reflection principle described above will be entailed by the reflection principle of $B_0(V_0)$. Thus we see from the foregoing considerations that the following is established.

\begin{theorem} Marshall's theory $B_0(V_0)$ proves that the height $\kappa$ of the universe $V_0$ is supercompact, and therefore by further application of the reflection principle that $V_0$ is a model for the assertion that there is a proper class of supercompact cardinals. \end{theorem}

\begin{proof} This is clear from the above discussion. \end{proof}

This resolves affirmatively the question asked by Marshall about whether $B_0(V_0)$ is sufficiently strong to prove the existence of a supercompact cardinal.

\bigskip

Later in the paper \cite{Marshall1989} Marshall introduces still stronger theories as motivation for the existence of $n$-huge cardinals. In the final section of her paper she goes even further up to the point of inconsistency with the axiom of choice. I hope in later work to elaborate on possible principled reasons, inspired by the lines of thought put forward by Marshall, and assuming that one does indeed regard the principles of Welch and Roberts as intrinsically justified, for accepting all large cardinals up to I1 and possibly still further beyond, while still finding principled reasons for not going to the point of inconsistency with choice.

\section{Concluding remarks}

Our exploration of the range of reflection principles which have been proposed in the literature have led us to outline two possible conceptions of intrinsic justification, one which carries us up roughly to the point of a remarkable cardinal, and is in line with Peter Koellner's conjecture that intrinsic justifications do not take us as far as an $\omega$-Erd\H{o}s cardinal, and one which extends up to the reflection principles of Welch and Roberts, and indeed, as suggested here, plausibly all the way up to a supercompact cardinal. We have offered an attempt to argue for philosophical theses (1), (2) and (3) as outlined in Section 1.

\pagebreak[4]


\begin{thebibliography}{99}

\bibitem{Boolos1971}
\newblock{George Boolos, 1971},
\newblock{The iterative conception of set},
\newblock{\em{Journal of Philosophy}},
\newblock{vol. 68, Issue 8, pp. 215--231.}

\bibitem{Cantor1932}
\newblock{Georg Cantor, 1932},
\newblock{\em{Gesammelte Abhandlungen mathematischen und philosophischen Inhalts}},
\newblock{Berlin: Springer, ed. E. Zermelo}

\bibitem{Cantor1883}
\newblock{Georg Cantor, 1883},
\newblock{\"Uber unendliche, lineare Punktmannigfaltigkeiten, 5, \em{Mathematische Annalen}},
\newblock{vol. 21, pp. 545--586. Appears in \cite{Cantor1932}.}

\bibitem{Chuaqui1978}
\newblock{R. Chuaqui, 1978},
\newblock{Bernays' class theory, ``Mathematical logic: proceedings of the first Brazilian conference"},
\newblock{(Campinas, 1977: A. I. Arruda et al., editors), Lecture Notes in Pure and Applied Mathematics},
\newblock{vol. 39, pp. 31--55, Marcel Dekker, New York.}

\bibitem{Drake1974}
\newblock{Frank Drake, 1974},
\newblock{Set Theory: An Introduction to Large Cardinals},
\newblock{North-Holland Publishing Company.}

\bibitem{Gitman2018}
\newblock{Victoria Gitman and Ralf Schindler, 2018},
\newblock{Virtual Large Cardinals, \em{Annals of Pure and Applied Logic}},
\newblock{vol. 69, Issue 12, pp. 1317--1334.}

\bibitem{Godel1990}
\newblock{Kurt G\"odel, 1990},
\newblock{\em{Collected Works, Volume II: Publications 1938-1974}},
\newblock{Oxford Univeristy Press, New York and Oxford, ed. Solomon Feferman, John W. Dawson, Jr., Stephen C. Kleene},
\newblock{Gregory H. Moore, Robert M. Solovay, and Jean van Heijenoort.}

\bibitem{Godel1964}
\newblock{Kurt G\"odel, 1964},
\newblock{What is Cantor's continuum problem?},
\newblock{In \cite{Godel1990}, pp. 254--270.}

\bibitem{Kanamori2003}
\newblock{Akihiro Kanamori, 2003},
\newblock{The Higher Infinite: Large cardinals in set theory from their beginnings, 2nd edition},
\newblock{Springer Monographs in Mathematics.}

\bibitem{Koellner2006}
\newblock{Peter Koellner, 2006}
\newblock{On the Question of Absolute Undecidability.}
\newblock{\em{Philosophia Mathematica},}
\newblock{vol. 14, number 2, pp. 153--188.}

\bibitem{Koellner2009}
\newblock{Peter Koellner, 2009},
\newblock{On Reflection Principles, \em{Annals of Pure and Applied Logic}},
\newblock{vol. 157, Issue 2--3, pp. 206--219.}

\bibitem{Marshall1989}
\newblock{M. Victoria Marshall R., 1989},
\newblock{Higher order reflection principles},
\newblock{\em{Journal of Symbolic Logic}},
\newblock{vol. 54, no. 2, pp. 474--489.}

\bibitem{McCallum2013}
\newblock{Rupert McCallum, 2013},
\newblock{A Consistency Proof for Some Restrictions of Tait's Reflection Principles},
\newblock{\em{Mathematical Logic Quarterly}},
\newblock{vol. 59, Issues 1--2, pp. 112--118.}

\bibitem{McCallum2020}
\newblock{Rupert McCallum, 2020},
\newblock{Extending ideas of Tait for incorporating higher-order parameters in schemes of reflection},
\newblock{\em{Journal of Applied Logics - IfCoLog Journal of Logics and their Applications}},
\newblock{vol. 7, no. 4, pp. 391--402.}

\bibitem{Potter2004}
\newblock{Michael Potter, 2004},
\newblock{Set Theory and Its Philosophy: A Critical Introduction},
\newblock{Oxford University Press.}

\bibitem{Reinhardt1974}
\newblock{William Reinhardt, 1974},
\newblock{Remarks on reflection principles, large cardinals, and elementary embeddings},
\newblock{Axiomatic Set Theory, Proc. Symp. Pure Math., ed. Dana Scott, vol. 13, pp. 189--205, A.M.S, Providence, RI.}

\bibitem{Roberts2017}
\newblock{Sam Roberts, 2017},
\newblock{A strong reflection principle},
\newblock{\em{The Review of Symbolic Logic}},
\newblock{vol. 10, Issue 4, pp. 651--662.}

\bibitem{Schindler1994}
\newblock{Ralf Schindler, 1994},
\newblock{A dilemma in the philosophy of set theory},
\newblock{\em{Notre Dame Journal of Formal Logic}},
\newblock{vol. 35, no. 3, pp. 458--463.}

\bibitem{Schindler2000}
\newblock{Ralf Schindler, 2000},
\newblock{Proper Forcing and Remarkable Cardinals},
\newblock{\em{The Bulletin of Symbolic Logic}},
\newblock{vol. 6, no. 2, pp. 176--184.}

\bibitem{Tait2005b}
\newblock {William Tait, 2005},
\newblock {\em{The Provenance of Pure Reason: Essays In the Philosophy of Mathematics and Its History}},
\newblock{Oxford University Press.}

\bibitem{Tait2005a}
\newblock {William Tait, 2005},
\newblock{Constructing Cardinals from Below. In \cite{Tait2005b}, Oxford University Press, pp. 133--154.}

\bibitem{Welch2017}
\newblock{Philip Welch, 2017},
\newblock{Obtaining Woodin's Cardinals},
\newblock{in ``Logic in Harvard: Conference celebrating the birthday of Hugh Woodin''},
\newblock{eds. A. Caicedo, J. Cummings, P. Koellner \& P. Larson, AMS Series, Contemporary Mathematics, vol. 690, pp. 161--176.}

\end{thebibliography}
\end{document}